\documentclass[12pt, oneside, emlines]{article}

\usepackage{tikz}
\usepackage{tkz-graph}
\usepackage{float}

\usepackage{orcidlink} 
\usepackage{amssymb,latexsym,xy,eucal,mathrsfs}
\usepackage[english]{babel}
\usepackage[utf8]{inputenc}
\usepackage{amsmath}
\usepackage{graphicx}
\usepackage[colorinlistoftodos]{todonotes}
\usepackage[affil-it]{authblk}
\usepackage{amssymb,latexsym,xy,eucal,mathrsfs,tikz}
\usepackage{amsthm,amssymb,amsfonts,amsmath}
\usepackage{latexsym,float,enumerate,graphics,graphicx}
\usepackage{hyperref}
\usepackage{multicol}
\usepackage{setspace}
\usepackage{mathtools}

\usepackage{indentfirst}
\usepackage{hyperref}
\usepackage{xcolor}
\hypersetup{
	colorlinks,
	linkcolor={blue!50!black},
	citecolor={green!50!black},
	urlcolor={blue!80!black}
}

\usepackage{mathtools}
\DeclarePairedDelimiter\ceil{\lceil}{\rceil}

\usepackage{geometry}

\theoremstyle{plain}
\newtheorem{lem}{Lemma}[section]
\newtheorem{lemma}[lem]{Lemma}
\newtheorem{theorem}[lem]{Theorem}

\newtheorem{prop}[lem]{Proposition}

\theoremstyle{definition}

\usepackage{lineno}

\newcommand{\x}{\times}

\usepackage{graphicx}
\usepackage{enumerate}

\usepackage[numbers,sort&compress]{natbib}
\begin{document}

\baselineskip 14truept

\title{The Domination and Secure Domination Numbers \\of Direct Product of Cliques with Paths and  Cycles}

\author[]{Deepak M. Bakal~\orcidlink{0009-0007-8505-759X} and S. A. Mane~\orcidlink{0000-0003-1214-918X}\\ { \tt {\normalsize iamdeepakbakal@gmail.com, manesmruti@yahoo.com}}}
\affil{\normalsize Department of Mathematics, Savitribai Phule Pune University, Pune 411007, India.}

\date{}

\footnotetext[1]{\textit{First Version 14/12/2023, submitted to journal for publication}}
\maketitle 

\begin{abstract}
 In this paper, we obtain the exact values of several domination parameters for the direct product of a complete graph with a path or a cycle. Specifically, we determine the domination number,  independent domination number,  $[1,2]$-domination number, secure domination number, and 2-domination number for this family of graphs. We show that, in these graphs, the independent domination number and the $[1,2]$-domination number coincide with the domination number, while the secure domination number coincides with the 2-domination number. Additionally, as a consequence of our findings, we provide counterexamples  to disprove some erroneous results in the literature.   
\end{abstract}

 	


\section{Introduction} 
\par Consider an undirected finite simple graph $G$ with vertex set $V(G)$ and edge set $E(G)$. A set $D \subseteq V(G)$ is called a \textit{dominating set} of $G$ if each vertex  $w \in V(G) \setminus D$ is adjacent to at least one vertex $v \in D$. In this case, we say that \textit{$v$ covers $w$}. We let $N(w)$ denote the \textit{open neighborhood} of $w$, that is, the set of vertices adjacent to $w$. The \textit{closed neighborhood} of $w$ is $N[w]=N(w) \cup \{w\}$. In other words, $D$ is a dominating set, if for each vertex  $w, N[w] \cap D$ is non-empty. If $A$ and $B$ are subsets of $V(G)$ we say that $A$ \textit{dominates} $B$ if every vertex of $B$ is adjacent to at least one vertex in  $A$ or is a vertex of $A$. The \textit{domination number} of a graph $G$ is the smallest cardinality of a dominating set and it is denoted by $\gamma(G)$. A dominating set represents a solution for numerous practical problems ranging from network science, to transportation streaming, to biological networks, to communication networks. Particularly, dominating sets are prominently used in the efficient organization of large-scale wireless ad hoc and sensor networks in the field of communication. Consider the following scenario:  suppose the vertices of a graph $G$ denote different locations and edges denote the link between them. We are interested in placing monitoring devices on some locations in such a way that all locations are under observation. Clearly it is desirable to do this with minimum number of devices. 

\par Cockayne \textit{et al.} \cite{cockayne2005protection} introduced the secure versions of dominating sets. A set $S \subseteq V$ is called a \textit{secure dominating set} of $G$ if it is a dominating set of $G$ and if, for each
$w \in V \setminus S$, there exists a vertex $v \in N(w) \cap S$ such that the \textit{swap} set $(S \setminus \{v\}) \cup \{w\}$ is again a
dominating set of $G$. In this case, we say that $v$ \textit{defends} $w$. The smallest cardinality of
a secure dominating set of $G$ is called the \textit{secure domination number} of $G$ and is denoted
by $\gamma_s(G)$. The name itself suggests that this concept is related to security. We illustrate this by extending the above example. Suppose that there is an attack at a location. If the location is placed on a device, then the device will resist the attack. Otherwise, we will move an adjacent device to resist the attack.  It is obviously desirable that the devices can still monitor all locations after moving. Such a problem can be modeled as the secure domination problem.
\par The notion of secure domination has been researched extensively. The decision version of computing the secure domination number is an NP-hard problem even when restricted to bipartite graphs, and split graphs \cite{merouane2015secure}, unit disc graphs \cite{wang2023algorithmic}. On the positive side there are linear-time algorithms for trees \cite{burger2014linear}, block graphs \cite{pradhan2018computing}, and cographs \cite{jha2019secure}, \cite{kivsek2021correcting}.  Cockayne \textit{et al.} \cite{cockayne2005protection} obtained exact values of $\gamma_s(G)$ for various graph classes, such as paths, cycles, complete
multipartite graphs. Rangel Hernández-Ortiz \textit{et al.} \cite{hernandez2021secure} studied secure domination number of the
rooted product graphs. Valveny and Rodriguez-Velazquez \cite{valveny2020protection} obtained the secure domination number
of the corona product of any graph with a discrete graph, the Cartesian product of two equal-sized
stars, and the Cartesian product of any clique other than $K_2$, with either a path, cycle, or star. This was extended by  Haythorpe \textit{et al.} \cite{haythorpe2022secure} by determing the secure domination numbers of the Cartesian products of $P_2, P_3$ with a path and a cycle, and for Mobius ladder graphs. 

\par A subset of vertices is \textit{independent}  if no two vertices in it are adjacent. An \textit{independent dominating set} of $G$ is a set that is
both dominating and independent in $G$. The \textit{independent domination number} $i(G)$ of $G$ is the minimum size of
an independent dominating set. It follows immediately that $\gamma(G) \leq i(G)$. A subset of vertices is a \textit{$[1,2]$-set} if every vertex $v \in V$ is either in it or  adjacent to at least one but no more than two vertices in it. The \textit{$[1,2]$-domination number} of a graph $G$ $\gamma_{_{[1, 2]}} (G)$ is the minimum cardinality of a $[1, 2]$-set of $G$. Obviously, any [1,2]-set is also a dominating set, so $\gamma(G) \leq \gamma_{_{[1, 2]}} (G)$.  A vertex subset of a graph $G$ is said to 2-dominate the graph if each vertex of $G$ is either in the given subset or  has
at least two neighbors in it. The minimum cardinality of a
2-dominating set, denoted by $\gamma_2(G)$,  is called the 2-domination number of the graph $G$.
A 2-dominating set is clearly a secure dominating set and so, $\gamma_2(G) \geq   \gamma_s(G)$.
\par Given two graphs $G$ and $H$, the \textit{direct product} $G \times H$  is the graph with vertex set $V (G) \times V (H)$ where two vertices $(g_1,h_1)$ and $(g_2,h_2)$ are adjacent
if and only if $g_1g_2\in E(G)$ and $h_1h_2 \in E(H)$. This product is also known as \textit{Kronecker product, tensor product, categorical product,  cardinal product, cross product, relational product} and \textit{graph conjunction}. The interest in direct product is due to the fact that
large networks such as the citation networks, neural networks, and internet graphs with several hundred million hosts, can be efficiently
modeled by the direct product \cite{leskovec2010kronecker}. For in-depth treatment of product graphs we refer the interested reader to book by Hammack \textit{et al.} \cite{hammack2011handbook}. Domination number of the direct product of graphs has been studied by a number of researchers. Brevsar \textit{et al.} \cite{brevsar2007dominating} obtained  an upper bound for the domination number of the direct product of graphs viz. $$\gamma(G \x H) \leq 3\gamma(G)\gamma(H).$$ Mekis \cite{mekivs2010lower} obtained a lower bound on the domination number viz.  $$  \gamma(G \x H) \geq \gamma(G)+ \gamma(H) - 1.  $$ Gravier \cite{gravier1995domination}  studied the domination number of direct product of a path with complement of another path. In this paper he presented a Vizing-like conjecture for direct product graphs, which was refuted by Klavzar \textit{et al.} \cite{klavẑar1996vizing}. Klobucar \cite{klobuvcar1999domination} \cite{klobuvcar1999domination2}, Cheriffi \textit{et al.} \cite{cherifi1999domination}  studied the domination number of direct product of two paths. Mekis \cite{mekivs2010lower}, Defant and Iyer \cite{defant2018domination}, Burcroff \cite{burcroff2018domination} and Vemuri \cite{vemuri2020domination} studied domination number of the direct product of finitely many cliques. In this paper, we focus on  domination and secure domination in the direct product of complete graphs with paths or cycles.
\par For a natural  number $n$, let $[n]$  denote the set $\{1, 2, . . . , n\}$. Let $P_n$ denote a path  with vertex set $[n]$, $C_n$ denote a cycle  with vertex set $[n]$ and $K_m$ denote a clique or a complete graph  with vertex set $[m]$. The graphs $P_n \x K_m$ and $C_n \x K_m$ have been of much interest. The total domination number of $C_n \x K_m$ and $P_n \x K_m$ was obtained by El-Zahar \cite{el2008total} and Zwierzchowski \cite{zwierzchowski2007total} respectively.
In \cite{lu2019identifying} Lu \textit{et al.} studied the identifying codes in $C_n \x K_m$ whereas in \cite{shinde2023identifying} Shinde \textit{et al.} studied the identifying codes in $P_n \x K_m$. Recently, Kuenzel and Rall \cite{kuenzel2023independent} gave the exact value of 
independent domination number of  $C_n \x K_m$ and $P_n \x K_m$. For undefined terminology and notation refer the book by West \cite{west2001introduction}.
\par   Note that, the notion of secure domination as presented in this paper differs from the concept of secure domination in graphs introduced by Brigham \textit{et al.} \cite{brigham2007security}. 

\subsection{Organization of the paper}
\par In this paper, we study the domination and secure domination in $C_n \x K_m$ and $P_n \x K_m$.
In Section 2, we discuss the notations and preliminaries. In Section 3, two basic lemmas are proved, which are crucial in the proof of our main results. In Section 4,  we produce minimal dominating sets of $C_n \x K_m$ and $P_n \x K_m$ and determine exact values of $\gamma(C_n \x K_m)$, $\gamma(P_n \x K_m)$. The dominating set constructed by us is an independent set as well as a [1,2]-set, so in turn, we obtain a large family of graphs with equal domination number, [1,2] domination number and independent domination number. 
The determination of the family of graphs $G$ for which $i(G) = \gamma(G)$ remains an open problem, as discussed in \cite{allan1978domination}.
Chellali et al. \cite{chellali20131} posed the question of identifying a family of graphs where $\gamma(G) = \gamma_{_{[1, 2]}} (G)$. 
We thus obtained a partial solution to these problems, distinct from the previously explored case of claw-free graphs \cite{allan1978domination}, \cite{chellali20131}. It is shown that, our results provide counterexamples to three erroneous claims in the literature about the domination numbers of direct product with path \cite{gravier1995domination}, \cite{Pnsitthiwirattham2012domination} and cycle \cite{Cnsitthiwirattham2012domination}. Finally, in Section 5 we study the secure domination in aforementioned family of graphs. We determine the exact values of $\gamma_s(C_n \x K_m)$, and $\gamma_s(P_n \x K_m)$. As a consequence, the exact value of 2-domination number of $C_n \x K_m$ and $P_n \x K_m$ is obtained.

\section{Preliminaries}
The vertex set of $P_n \x K_m$ and $C_n \x K_m$ is $[n] \x [m]$ and henceforth will be denoted by $V$. It can be partitioned into disjoint arrays of vertices, hereafter referred to as \textit{columns} and \textit{rows}. A column and a row of $V$ 
are all the vertices with the same first and second coordinate, respectively. 
For $i \in [n]$, the $i^{th}$ \textit{column} is $X_i = \{(i, j) \mid~ j \in [m]\}.$ Similarly for $j \in [m]$, the $j^{th}$ \textit{row} is $R_j = \{(i, j) \mid~ i \in [n]\}.$ For both graph $P_n \x K_m$, and $C_n \x K_m$ rows are presented horizontally and columns are presented vertically. The following figure shows the graph $P_6 \x K_5$ and $C_6 \x K_5$. Graph $P_6 \x K_5$ is a subgraph of $C_6 \x K_5$. The additional edges of $C_6 \x K_5$ are shown by dotted line.

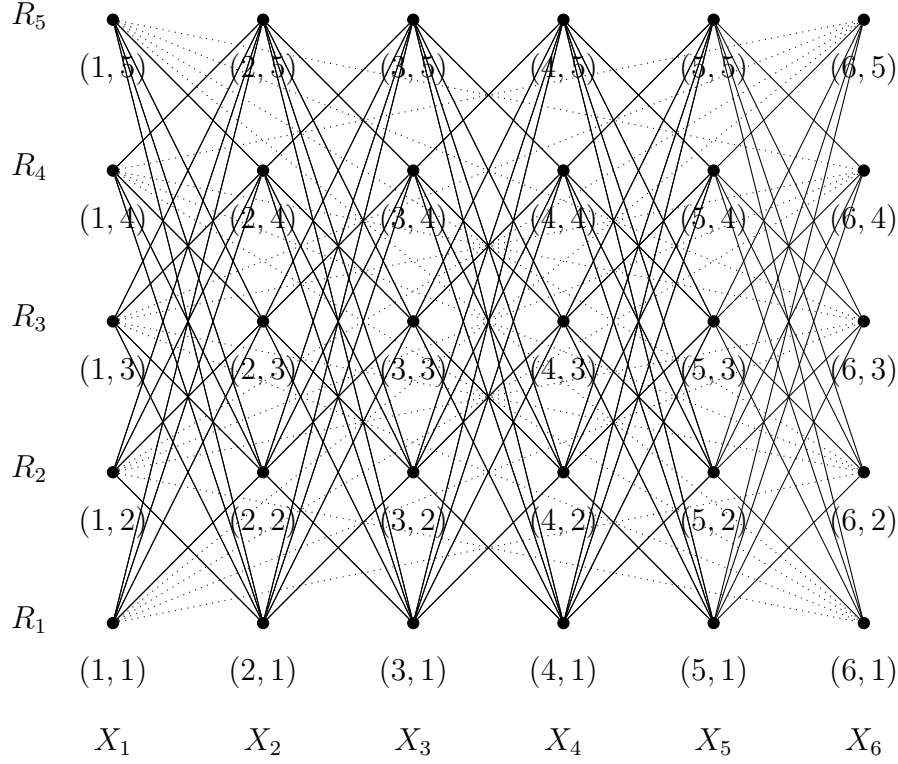
\begin{figure}
\begin{center}
	\noindent
	\begin{tikzpicture}

		\foreach \i in {0,...,5} {
			\node  at (2*\i, 1) {};
			\pgfmathtruncatemacro{\j}{\i+1}
			\node[below, yshift=0.0mm] at (2*\i, 0.75) {$X_{\j}$};
			
		}

		\foreach \i in {1,...,5} {
			\node (P1\i) at (0, 2*\i) {};
			\node[right, yshift=0.5mm] at (-1.5, 2*\i) {$R_\i$};
			
		}
		
		\foreach \i in {1,...,5} {
			\node[circle, draw, fill=black, inner sep=1.5pt, minimum size=2pt] (P1\i) at (0, 2*\i) {};
			\node[below, yshift=-3mm] at (0, 2*\i) {$(1, \i)$};
		}
		
		\foreach \i in {1,...,5} {
			\node[circle, draw, fill=black, inner sep=1.5pt, minimum size=2pt] (P2\i) at (2, 2*\i) {};
			\node[below, yshift=-3mm] at (2, 2*\i) {$(2, \i)$};
		}
		
		\foreach \i in {1,...,5} {
			\node[circle, draw, fill=black, inner sep=1.5pt, minimum size=2pt] (P3\i) at (4, 2*\i) {};
			\node[below, yshift=-3mm] at (4, 2*\i) {$(3, \i)$};
		}

		\foreach \i in {1,...,5} {
			\node[circle, draw, fill=black, inner sep=1.5pt, minimum size=2pt] (P4\i) at (6, 2*\i) {};
			\node[below, yshift=-3mm] at (6, 2*\i) {$(4, \i)$};
		}
		
		\foreach \i in {1,...,5} {
			\node[circle, draw, fill=black, inner sep=1.5pt, minimum size=2pt] (P5\i) at (8, 2*\i) {};
			\node[below, yshift=-3mm] at (8, 2*\i) {$(5, \i)$};
		}
		
		\foreach \i in {1,...,5} {
			\node[circle, draw, fill=black, inner sep=1.5pt, minimum size=2pt] (P6\i) at (10, 2*\i) {};
			\node[below, yshift=-3mm] at (10, 2*\i) {$(6, \i)$};
		}

		\foreach \i in {1,...,5} {
			\foreach \j in {1,...,5} {
				\ifnum\i>\j
				\draw (P1\i) -- (P2\j);
				\draw[dotted] (P1\i) -- (P6\j);
				\fi
			}
		}

		\foreach \i in {1,...,5} {
			\foreach \j in {1,...,5} {
				\ifnum\i<\j
				\draw (P1\i) -- (P2\j);
				\draw[dotted] (P1\i) -- (P6\j);
				\fi
			}
		}

		
		\foreach \i in {1,...,5} {
			\foreach \j in {1,...,5} {
				\ifnum\i>\j
				\draw (P2\i) -- (P1\j);
				\draw (P2\i) -- (P3\j);
				\fi
			}
		}

		\foreach \i in {1,...,5} {
			\foreach \j in {1,...,5} {
				\ifnum\i<\j
				\draw (P2\i) -- (P1\j);
				\draw (P2\i) -- (P3\j);
				\fi
			}
		}

		
		\foreach \i in {1,...,5} {
			\foreach \j in {1,...,5} {
				\ifnum\i>\j
				\draw (P3\i) -- (P2\j);
				\draw (P3\i) -- (P4\j);
				\fi
			}
		}

		\foreach \i in {1,...,5} {
			\foreach \j in {1,...,5} {
				\ifnum\i<\j
				\draw (P3\i) -- (P2\j);
				\draw (P3\i) -- (P4\j);
				\fi
			}
		}

		
		\foreach \i in {1,...,5} {
			\foreach \j in {1,...,5} {
				\ifnum\i>\j
				\draw (P4\i) -- (P3\j);
				\draw (P4\i) -- (P5\j);
				\fi
			}
		}

		\foreach \i in {1,...,5} {
			\foreach \j in {1,...,5} {
				\ifnum\i<\j
				\draw (P4\i) -- (P3\j);
				\draw (P4\i) -- (P5\j);
				\fi
			}
		}

		\foreach \i in {1,...,5} {
			\foreach \j in {1,...,5} {
				\ifnum\i>\j
				\draw (P5\i) -- (P4\j);
				\draw (P5\i) -- (P6\j);
				\fi
			}
		}

		\foreach \i in {1,...,5} {
			\foreach \j in {1,...,5} {
				\ifnum\i<\j
				\draw (P5\i) -- (P4\j);
				\draw (P5\i) -- (P6\j);
				\fi
			}
		}

	\end{tikzpicture}
	
	\label{fig:myfigurefirst}
         \caption{Graph $P_6 \times K_5$ and $C_6 \times K_5$}
	
\end{center}
\end{figure}

\par It is easy to observe that any two vertices of $P_n  \times K_m$ are adjacent if and only if they belong to different rows and consecutive columns. Also, $C_n  \times K_m$ is a supergraph of $P_n \x K_m$ with same vertex set and additional edges 
 between first and last columns, $X_1$ and $X_n$, respectively.



\par \noindent The following proposition was proved by Cockayne \textit{et al.} \cite{cockayne2005protection}.
\begin{prop}[\cite{cockayne2005protection} Theorem 12] \label{Prop2}
a. For $n \geq 1$,	$\gamma_s(P_n) = \ceil*{\frac{3n}{7}} $. \\
	
b. For $n \geq 4$,	$\gamma_s(C_n) = \ceil*{\frac{3n}{7}} $. \\

c. $\gamma_s(C_2) = \gamma_s(C_3) = 1$.

\end{prop}

The path, cycle and clique on two vertices are isomorphic.  The direct product of a graph $G$ with $K_2$ is called  the \textit{bipartite double cover of} $G$. We obtain the secure domination number of bipartite double cover of $P_n$ and $C_n$.

\begin{prop} For $n \geq 2$, 
	
	\begin{enumerate}
		\item $\gamma( P_n \x K_2) = 2    \ceil*{\frac{n}{3}}$.
		\item $\gamma( C_n \x K_2) =
		\begin{cases} 
			\ceil*{\frac{2n}{3}}, &\mbox{\textnormal{\quad if}}~ n \textnormal {~is~odd};\\
			2\ceil*{\frac{n}{3}}, &\mbox{\textnormal{\quad if}}~ n  \textnormal {~is even}.
		\end{cases}$

		\item $\gamma_s( P_n \x K_2) = 2    \ceil*{\frac{3n}{7}}$. 
		\item $\gamma_s( C_n \x K_2) =
		\begin{cases} 
			\ceil*{\frac{6n}{7}}, &\mbox{\textnormal{\quad if}}~ n \textnormal {~is~odd};\\
			2\ceil*{\frac{3n}{7}}, &\mbox{\textnormal{\quad if}}~ n  \textnormal {~is even}.
		\end{cases}$

	\end{enumerate}
\end{prop}
\begin{proof}
	The graph $P_n \x K_2$ is isomorphic to two disjoint copies of $P_n$. If $n$ is even, the graph $ C_n \x K_2 $ is isomorphic to two disjoint copies of $C_n$. If  $n$ is odd, the graph $ C_n \x K_2 $ is isomorphic to $C_{2n}$. Hence the proposition follows from   Proposition \ref{Prop2} and the fact $\gamma(P_n) = \gamma(C_n) = \ceil*{\frac{n}{3}}$.
\end{proof}


\section{Some Lemmas}

\par   In this section, we prove some essential lemmas which will be used in the proof of our main results appearing in the subsequent sections.
 For a dominating set $D$ of either $P_n \x K_m$ or  $C_n \x K_m$, and for $i \in [n]$, let $d_i$ denote the number of dominating vertice in $X_i$ \textit{i.e.} the cardinality of $D \cap X_i$. In $P_n \x K_m$, it is convenient to assume that $X_0=X_{n+1}=\emptyset$ (\emph{i.e.} $d_0=d_{n+1}=0$), whereas in case of $C_n \x K_m$, we suppose $X_n=X_{0}$ and $X_{n+1}=X_{1}$ (\emph{i.e.} $d_n=d_{0}$ and $d_{n+1}=d_{1}$).\\



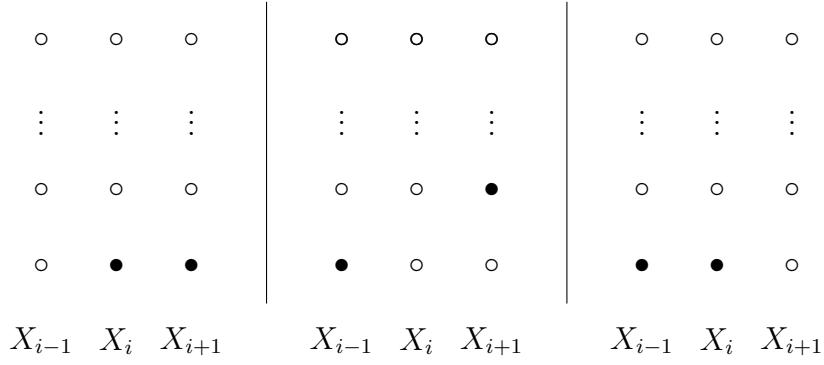
\begin{figure}
\begin{center}
\noindent
	
	\begin{tikzpicture}
		
		\node at (0,3) {$\circ$};
		\node at (1,3) {$\circ$};
		\node at (2,3) {$\circ$};

		\node at (-2,3) {$\circ$};
		\node at (-3,3) {$\circ$};
		\node at (-4,3) {$\circ$};

		\node at (-2,2) {$\vdots$};
		\node at (-3,2) {$\vdots$};
		\node at (-4,2) {$\vdots$};
		
		\node at (-2,1) {$\circ$};
		\node at (-3,1) {$\circ$};
		\node at (-4,1) {$\circ$};
		
		\node at (-2,0) {$\bullet$};
		\node at (-3,0) {$\bullet$};
		\node at (-4,0) {$\circ$};
		
		\node at (-2,-1)  {$X_{i+1}$};
		\node at (-3,-1)  {$X_{i}$};
		\node at (-4,-1)  {$X_{i-1}$};
		
		\draw (-1,-0.5) -- (-1,3.5);
		
		\node at (0,3) {$\circ$};
		\node at (1,3) {$\circ$};
		\node at (2,3) {$\circ$};

		\node at (0,2) {$\vdots$};
		\node at (1,2) {$\vdots$};
		\node at (2,2) {$\vdots$};

		\node at (0,1) {$\circ$};
		\node at (1,1) {$\circ$};
		\node at (2,1) {$\bullet$};
		
		\node at (0,0) {$\bullet$};
		\node at (1,0) {$\circ$};
		\node at (2,0) {$\circ$};
		
		\node at (0,-1) {$X_{i-1}$};
		\node at (1,-1) {$X_{i}$};
		\node at (2,-1) {$X_{i+1}$};
		
		\draw (3,-0.5) -- (3,3.5);

		\node at (4,3) {$\circ$};
		\node at (5,3) {$\circ$};
		\node at (6,3) {$\circ$};

		\node at (4,2) {$\vdots$};
		\node at (5,2) {$\vdots$};
		\node at (6,2) {$\vdots$};

		\node at (4,1) {$\circ$};
		\node at (5,1) {$\circ$};
		\node at (6,1) {$\circ$};
		
		\node at (4,0) {$\bullet$};
		\node at (5,0) {$\bullet$};
		\node at (6,0) {$\circ$};
		
		\node at (4,-1) {$X_{i-1}$};
		\node at (5,-1) {$X_{i}$};
		\node at (6,-1) {$X_{i+1}$};
		
	\end{tikzpicture}
\caption{Dominating Sets of $X_i$}
				\label{fig:myfigure}
			\end{center}
\end{figure}

\begin{lemma} \label{lemma2}
	Let $m \geq 3$ and $n \geq 2$. Then for a dominating set $D$ of $P_n \x K_m$ or  $C_n \x K_m$,   $d_{i-1}+d_{i}+d_{i+1} \geq 2,$ for $1 \leq i \leq n.$
\end{lemma}

\begin{proof}
	Let  $A$ be the subset of dominating set $D$ which dominates the column $X_i$. Since $A$ dominates $X_i$, the vertex $(i,1)$ is either in $A$ or there exists a vertex, say $y$, 
	that covers $(i,1)$. Let vertex $y$ lies in the row $R_j$. Now, the vertex $(i,j)$ is not dominated by either $y$ or $(1,j)$. So, there exists another vertex in $A$ which dominates $(i,j)$. Hence $|A| \geq 2$.  Now, the column $X_i$ is dominated only by the vertices in $(X_{i-1} \cup  X_i \cup X_{i+1})$. So $A = D \cap (X_{i-1} \cup X_i \cup X_{i+1})$. It is easy to see that, $|A|= |D \cap X_{i-1}| + |D \cap X_{i}| + |D \cap X_{i+1}| = d_{i-1}+d_{i}+d_{i+1}$. Hence, we obtain, $d_{i-1}+d_{i}+d_{i+1} \geq 2,$  for $1 \leq i \leq n.$
\end{proof}


\par  In the next lemma, we identify doubleton sets 
which dominate column $X_i$.
\begin{lemma} \label{SLemma1}
	Let  $A \subseteq V$, say   $A = \{a,b\}$. Then $A$ dominates $X_i$ if and only if 
	
	\begin{enumerate}
		\item $A \subseteq X_{i-1} \cup X_{i} \cup X_{i+1}$ with $A \nsubseteq X_i$; and
		\item Vertices $a,b$ lie in consecutive columns if and only if both of them lie in the same row.
	\end{enumerate}

\end{lemma}  

\begin{proof}
($\Rightarrow$)  \begin{enumerate}
	
	\item Vertices in $X_i$ are adjacent to vertices in $X_{i-1}$ and $X_{i+1}$ alone. As $A$ dominates $X_i$, clearly $A \subseteq X_{i-1} \cup X_{i} \cup X_{i+1}$ by the adjacency relation. Since $|X_i|= m \geq 3$, obviously $A \nsubseteq X_i$; because otherwise $A$ does not dominate vertices in $X_i \setminus \{a,b\}$.

	\item Suppose $a,b$ lie in consecutive columns $X_{i-1}, X_{i}$. Say $a = (i-1, p),~~ b =(i,q) $. If $p \neq q$, then the vertex $(i,p)$ is not dominated, so $p=q$ must hold, \emph{i.e.} $a,b$ lie in same row. The other case $a,b$ lying in consecutive columns $X_{i}, X_{i+1}$ is similar.

	\par Suppose $a,b$  do not lie in consecutive columns. Then,  either (a) they lie in same column or (b) one vertex lies in $X_{i-1}$ and other vertex $X_{i+1}.$ We consider cases (a) and (b) as follows.\begin{enumerate}[(a)]
		\item If $A \subseteq X_{i-1}$ or $A \subseteq X_{i+1}$, then obviously $a,b$ do not lie in the same row.
		\item Without loss of generality, suppose $a \in X_{i-1}$ and $b \in X_{i+1}$. If $a=(i-1,p)$ and $b = (i+1, q)$ with $p=q$ then the vertex $(i,p) \in X_i$ is not dominated.
	\end{enumerate}

\end{enumerate}	 
	
($\Leftarrow$) This is easy to observe (see Figure \ref{fig:myfigure}).

\end{proof}

\section{Domination Number of  \boldmath  $C_n \x K_m$ and $P_n \x K_m$ }

\subsection{Domination Number of  \boldmath  $C_n \x K_m$}
In this section, we determine the domination number of the graph  $C_n \x K_m$. We do this by obtaining an upper bound on $\gamma(C_n \x K_m)$, by producing a dominating set. We subsequently prove that the derived upper bound is, in fact, the domination number $\gamma(C_n \x K_m)$ by establishing it as a lower bound on the size of dominating set  of $C_n\times K_m$.\\ 

\par \noindent 
The next proposition establishes an upper bound on $\gamma(C_n\times K_m)$, for $n \geq 6$. The cases $n \leq 5$ are discussed at the end of this section.

\begin{prop}\label{CLemma2} For $m \geq 3,~n\geq 6$, 
	$$\gamma(C_n\times K_m)\leq \begin{cases} 
		2k, &\mbox{\textnormal{\quad if}}~ n=3k;\\
		2k+1, &\mbox{\textnormal{\quad if}}~ n=3k+1;\\
		2k+2, &\mbox{\textnormal{\quad if}}~ n=3k+2.
	\end{cases}$$ 
\end{prop} 
\begin{proof} 
	Consider the following cases. In each case, we produce a dominating set of the required size by partitioning the vertex set, as shown in the following figure \ref{fig:DomCnKm}.
	\begin{figure}[H]

		\begin{center}
				\noindent 
				\begin{tikzpicture}[scale=1]

		\node at (0,1) {$\tiny R_{1}$};
		\node at (0,2) {$R_{2}$};
		\node at (0,3) {$R_{3}$};
		\node at (0,4) {$\vdots$};
		\node at (0,5) {$R_{m}$};

		\node at (1,0) {$X_{1}$};
		\node at (2,0) {$X_{2}$};
		\node at (3,0) {$X_{3}$};
		
		\draw (3.5,0.5) -- (3.5,5.5);
		
		\node at (4,0) {$X_{4}$};
		\node at (5,0) {$X_{5}$};
		\node at (6,0) {$X_{6}$};
		
		\draw (6.5,0.5) -- (6.5,5.5);
		
		\node at (7,1) {$\ldots$};
		\node at (7,5) {$\ldots$};
		\node at (7,4) {$\ddots$};
		
		\draw (7.5,0.5) -- (7.5,5.5);

		\node at (8,0) {$X_{_{3k-2}}$};
		\node at (9,0) {$X_{_{3k-1}}$};
		\node at (10,0) {$X_{_{3k}}$};
		
		\draw (10.5,0.5) -- (10.5,5.5);
		
		\node at (11,0) {$X_{_{3k+1}}$};
		
		\draw (11.5,0.5) -- (11.5,5.5);
		
		\node at (12,0) {$X_{_{3k+2}}$};
		
		\node at (1,1) {$\bullet$};
		\node at (1,2) {$\circ$};
		\node at (1,3) {$\circ$};
		\node at (1,4) {$\vdots$};
		\node at (1,5) {$\circ$};
		
		\node at (2,1) {$\bullet$};
		\node at (2,2) {$\circ$};
		\node at (2,3) {$\circ$};
		\node at (2,4) {$\vdots$};
		\node at (2,5) {$\circ$};
		
		\node at (3,1) {$\circ$};
		\node at (3,2) {$\circ$};
		\node at (3,3) {$\circ$};
		\node at (3,4) {$\vdots$};
		\node at (3,5) {$\circ$};

		\node at (4,1) {$\circ$};
		\node at (4,2) {$\bullet$};
		\node at (4,3) {$\circ$};
		\node at (4,4) {$\vdots$};
		\node at (4,5) {$\circ$};
		
		\node at (5,1) {$\circ$};
		\node at (5,2) {$\bullet$};
		\node at (5,3) {$\circ$};
		\node at (5,4) {$\vdots$};
		\node at (5,5) {$\circ$};
		
		\node at (6,1) {$\circ$};
		\node at (6,2) {$\circ$};
		\node at (6,3) {$\circ$};
		\node at (6,4) {$\vdots$};
		\node at (6,5) {$\circ$};

		\node at (8,1) {$\circ$};
		\node at (8,2) {$\circ$};
		\node at (8,3) {$\bullet$};
		\node at (8,4) {$\vdots$};
		\node at (8,5) {$\circ$};
		
		\node at (9,1) {$\circ$};
		\node at (9,2) {$\circ$};
		\node at (9,3) {$\bullet$};
		\node at (9,4) {$\vdots$};
		\node at (9,5) {$\circ$};
		
		\node at (10,1) {$\circ$};
		\node at (10,2) {$\circ$};
		\node at (10,3) {$\circ$};
		\node at (10,4) {$\vdots$};
		\node at (10,5) {$\circ$};

		\node at (11,1) {$\bullet$};
		\node at (11,2) {$\circ$};
		\node at (11,3) {$\circ$};
		\node at (11,4) {$\vdots$};
		\node at (11,5) {$\circ$};

		\node at (12,1) {$\bullet$};
		\node at (12,2) {$\circ$};
		\node at (12,3) {$\circ$};
		\node at (12,4) {$\vdots$};
		\node at (12,5) {$\circ$};

	\end{tikzpicture}	
                    \caption{Dominating Sets of $C_n \times K_m$}

		\end{center}
	
	\label{fig:DomCnKm}
	\end{figure}
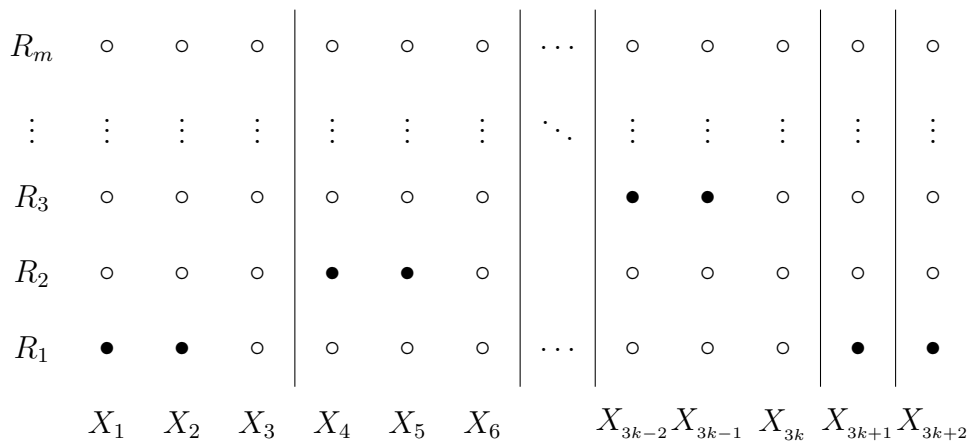

	\begin{itemize} 
		\item Case 1. Let $n=3k$.\\
		In this case we partition the vertex set $V$ as follows. Let $B_i$ denote the block of three columns $\left(X_{3i-2} \cup\nolinebreak X_{3i-1} \cup\nolinebreak X_{3i}\right)$.
		
		Write, $V =  \bigcup\limits_{i=1}^{k} (X_{3i-2} \cup X_{3i-1} \cup X_{3i})	 = \bigcup\limits_{i=1}^{k} B_i.$\\
		
		From each $B_i$, with $1 \leq k-1$, we choose two vertices $(3i-2,j), (3i-1,j)$ where \[
		j = \begin{cases}
			1 & \text{if $i$ is odd} \\
			2 & \text{if $i$ is even}
		\end{cases}.
		\] and from $B_k$, we choose two vertices $(3k-2,3), (3k-1,3)$, to get the set, \[D_{3k} = \left\{(1, 1), (2, 1), (4, 2), (5, 2), \ldots, (3k-2, 3), (3k-1, 3)\right\}.\]
		Observe that $D_{3k}$ is a dominating set of $C_{3k} \x K_m. $
		Hence,  $\gamma(C_{3k} \x K_m) \leq 2k.$

		\item Case 2. Let $n=3k+1$.\\
		Consider the set $D_{3k+1} = D_{3k} \cup \{(3k+1, 1)\}$.\\
		Observe that $D_{3k+1}$ is a dominating set of $C_{3k+1} \x K_m$.\\
      	Hence,  $\gamma(C_{3k+1} \x K_m) \leq 2k+1.$ 
		
		\item Case 3. Let $n=3k+2$.\\
		Consider the set $D_{3k+2} = D_{3k+1} \cup \{(3k+2,1)\}=D_{3k} \cup \{(3k+1,1),(3k+2,1)\}$.\\
		Observe that $D_{3k+2}$ is a dominating set of $C_{3k+2} \x K_m.$\\
		Hence,  $\gamma(C_{3k+2} \x K_m) \leq 2k+2.$ 				\end{itemize}
	Therefore, the proposition is proved.\end{proof}

We now estalish a lower bound on $\gamma(C_n\times K_m)$ by partitioning the vertex set into blocks consisting of three consecutive columns and invoking inequality in Lemma \ref{lemma2}.

\begin{theorem}\label{CLemma1} For $m \geq 3,~n\geq 6$, 
	$$\gamma(C_n\times K_m) = \begin{cases} 
		2k, &\mbox{\textnormal{\quad if}}~ n=3k;\\
		2k+1, &\mbox{\textnormal{\quad if}}~ n=3k+1;\\
		2k+2, &\mbox{\textnormal{\quad if}}~ n=3k+2.
	\end{cases}$$ 
\end{theorem} 

\begin{proof}
	Let $m \geq 3, n \geq 2$ and $D$ be a dominating set of $C_n  \times K_m$.
	
	\begin{itemize}
		\item Case 1. Let $n=3k$.
		$$|D| = \sum_{i=1}^{n} d_i = \sum_{i=1}^{k} (d_{3i-2}+d_{3i-1}+d_{3i}).$$
		Now, by Lemma \ref{lemma2}, $(d_{3i-2}+d_{3i-1}+d_{3i}) \geq 2$, for $1 \leq i \leq k$.
		Hence $|D| \geq 2k$.

		\item Case 2. Let $n=3k+1$.\\
		We prove that, if $ \left| D \right| \leq 2k$ then $D$ is not a dominating set. On the contrary, suppose $D$ is a dominating set such that $\left| D \right| \leq 2k$.
		
		Excluding the column $X_j$, partition the remaining vertex set $V$ into $k$ number of blocks of three consecutive columns. So by Lemma \ref{lemma2}, we get,
		
		$$\left( \sum_{i=1}^{3k+1} d_i \right) - d_j \geq 2k.$$ 
		
		Now we have assumed that, $$|D| = \left( \sum_{i=1}^{3k+1} d_i \right) \leq 2k,$$ we get $d_j = 0$. Since $j$ is arbitrary, we arrive at a contradiction. \\Hence, $|D| \geq 2k+1$.
		
%
%
%
%
		
%

		\item Case 3. Let $n=3k+2$.\\
		We will prove that any set with cardinality less than or equal to $2k+1$ is not a dominating set. 
		On the contrary, suppose $D$ is a dominating set such that $\left| D \right| \leq 2k+1$.
		
		Excluding the two consecutive columns $X_j$ and $X_{j+1}$, partition the remaining vertex set $V$ into $k$ number of blocks of three consecutive columns. So by Lemma \ref{lemma2}, we get,
		
		$$\left( \sum_{i=1}^{3k+1} d_i \right) - (d_j+d_{j+1}) \geq 2k.$$ 
		
		Now we have assumed that, $$|D| = \left( \sum_{i=1}^{3k+1} d_i \right) \leq 2k+1,$$ we obtain, $$d_j+d_{j+1} \leq 1.$$
		Similarly,  $$d_{j-1}+d_j \leq 1.$$ Therefore, $$d_{j-1}+d_j+d_{j+1}+d_j \leq 2$$
		
	But by Lemma \ref{lemma2} $$d_{j-1}+d_j+d_{j+1} \geq 2.$$
	So, we get that $d_j = 0$. Since $j$ is arbitrary, we arrive at a contradiction.\\ Hence, $|D| \geq 2k+2$.
	\end{itemize}  
	\noindent Therefore, for $m \geq 3$ and $n \geq 2$ $$\gamma(C_n\times K_m) \geq \begin{cases} 
	2k, &\mbox{\textnormal{\quad if}}~ n=3k;\\
	2k+1, &\mbox{\textnormal{\quad if}}~ n=3k+1;\\
	2k+2, &\mbox{\textnormal{\quad if}}~ n=3k+2.
	\end{cases}$$ 

\noindent Now if we assume $n \geq 6$, the theorem follows using Propositon \ref{CLemma2}. \end{proof}


\noindent In the next proposition we determine the value of $C_n \x K_m$, for $n < 6 $.

\begin{prop}\label{prop2345}

\noindent For $m \geq 3$ $$\gamma(C_n\times K_m) = \begin{cases} 
													2, &\mbox{\textnormal{\quad if}}~ n=2 ~~~~~~ 
													\\
													3, &\mbox{\textnormal{\quad if}}~ n=3 ~~~~~~ 
													\\
													4, &\mbox{\textnormal{\quad if}}~ n=4,5.
\end{cases}$$ 	
\end{prop}

\begin{proof}
The cycle graphs $C_2$ and $C_3$ are isomorphic to the complete graphs $K_2$ and $K_3$ respectively. The domination number $\gamma(K_2 \x K_m)$ and $\gamma(K_3 \x K_m)$ are  obtained by Mekis [Proposition 2.3, \cite{mekivs2010lower}].
\par Consider $C_4 \x K_m$, the proof of Theorem \ref{CLemma1} gives that $\gamma(C_4 \x K_m) \geq 3$. We prove that set of order three can not be a dominating set of $C_4 \x K_m$. Let, if possible, there exists a dominating set $D$ such that $|D|=d_1+d_2+d_3+d+4 = 3$. Without loss, we assume that $d_4=0$. 
Now, by Lemma \ref*{lemma2}, $d_4+d_2+d_3=d_2+d_3\geq 2$. Similarly $d_3+d_1\geq 2$ and $d_1+d_2\geq 2$. But $d_1+d_2+d_3=3$, so we get $d_1 = d_2 = d_3 =1$. Let $D = \{(1,i), (2,j), (3,k)\}$. If all the three $i,j,k$ are equal , then the vertex $(4,i)$ is not dominated. If all the three $i,j,k$ are distinct, then the vertices $(1,j), (3,j)$ are not dominated. The last possibility is exacty two of $i,j,k$ are equal, say $i=j$ and $i \neq k$, then $(3,i)$ is not dominated. Thus we get a contradiction. Therefore $\gamma(C_4\times K_m) \geq 4$. It is easy to see that the first row $R_1$ dominates $C_4 \x K_m$. Hence $\gamma(C_4 \x K_m) = 4$.
\par Consider $C_5 \x K_m$, the previous proof implies $\gamma(C_5 \x K_m) \geq 4$. It is easy to see that the set $\{(1,1), (2,1), (4,1), (5,1)\}$
is a dominating set of $C_5 \x K_m$. Hence $\gamma(C_5 \x K_m) = 4$.\end{proof}

%

\subsection{Domination Number of  \boldmath  $P_n \x K_m$}

In this section, we determine the domination number of the graph  $P_n \x K_m$. The idea is similar to the one discussed in preceding section. The graph $P_n  \times K_m$ is a subgraph of $C_n \x K_m$ with vertex set $V $ except edges drawn between first and last columns, $X_1$ and $X_n$, respectively. As noted previously, we suppose $d_0=d_{n+1}=0$. As path and cycle on two vertices are isomorphic, $\gamma(P_2 \x K_m)$ is obtained in Proposition \ref{prop2345}. So we can assume $n \geq 3$.\\

\noindent First we establish an upper bound on $\gamma(P_n\times K_m)$.

\begin{prop}\label{PLemma2} 
	For $m \geq 3,~n \geq 3$, 
	$$\gamma(P_n\times K_m)\\ \leq \begin{cases} 
		2k+1 &\mbox{\textnormal{if}}~ n=3k,\\
		2k+2 &\mbox{\textnormal{if}}~ n=3k+1,~ 3k+2.
	\end{cases}$$ 
\end{prop}	
\begin{proof}
	Consider the following cases. In each case, we produce a dominating set of the required size.
	\begin{itemize} 
		\item Case 1. Let $n=3k+1$.\\
		In this case we partition $V$ as follows, shown in the figure \ref{DomP3kplus1Km} .
		
		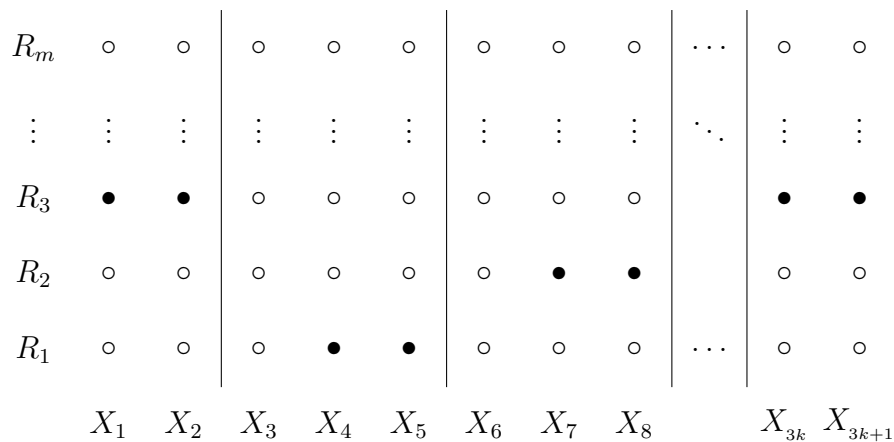
\begin{figure}[H]

		\begin{center}
			\noindent 
			\begin{tikzpicture}[scale=1]

		\node at (0,1) {$\tiny R_{1}$};
		\node at (0,2) {$R_{2}$};
		\node at (0,3) {$R_{3}$};
		\node at (0,4) {$\vdots$};
		\node at (0,5) {$R_{m}$};

		\node at (1,0) {$X_{1}$};
		\node at (2,0) {$X_{2}$};
		\node at (3,0) {$X_{3}$};
		
		\draw (2.5,0.5) -- (2.5,5.5);
		
		\node at (4,0) {$X_{4}$};
		\node at (5,0) {$X_{5}$};
		\node at (6,0) {$X_{6}$};
		
		\node at (7,0) {$X_{7}$};
		\node at (8,0) {$X_{8}$};

		\draw (5.5,0.5) -- (5.5,5.5);
		
		
		\draw (8.5,0.5) -- (8.5,5.5);

		\node at (10,0) {$X_{_{3k}}$};
		
		
		\node at (11,0) {$X_{_{3k+1}}$};
		
%

		\node at (1,1) {$\circ$};
		\node at (1,2) {$\circ$};
		\node at (1,3) {$\bullet$};
		\node at (1,4) {$\vdots$};
		\node at (1,5) {$\circ$};
		
		\node at (2,1) {$\circ$};
		\node at (2,2) {$\circ$};
		\node at (2,3) {$\bullet$};
		\node at (2,4) {$\vdots$};
		\node at (2,5) {$\circ$};
		
		\node at (3,1) {$\circ$};
		\node at (3,2) {$\circ$};
		\node at (3,3) {$\circ$};
		\node at (3,4) {$\vdots$};
		\node at (3,5) {$\circ$};

		\node at (4,1) {$\bullet$};
		\node at (4,2) {$\circ$};
		\node at (4,3) {$\circ$};
		\node at (4,4) {$\vdots$};
		\node at (4,5) {$\circ$};
		
		\node at (5,1) {$\bullet$};
		\node at (5,2) {$\circ$};
		\node at (5,3) {$\circ$};
		\node at (5,4) {$\vdots$};
		\node at (5,5) {$\circ$};

		\node at (6,1) {$\circ$};
		\node at (6,2) {$\circ$};
		\node at (6,3) {$\circ$};
		\node at (6,4) {$\vdots$};
		\node at (6,5) {$\circ$};

		\node at (7,1) {$\circ$};
		\node at (7,2) {$\bullet$};
		\node at (7,3) {$\circ$};
		\node at (7,4) {$\vdots$};
		\node at (7,5) {$\circ$};
		
		\node at (8,1) {$\circ$};
		\node at (8,2) {$\bullet$};
		\node at (8,3) {$\circ$};
		\node at (8,4) {$\vdots$};
		\node at (8,5) {$\circ$};
		
		

          \node at (9,1) {$\ldots$};
          \node at (9,5) {$\ldots$};
          \node at (9,4) {$\ddots$};

	\draw (9.5,0.5) -- (9.5,5.5);

		\node at (10,1) {$\circ$};
		\node at (10,2) {$\circ$};
		\node at (10,3) {$\bullet$};
		\node at (10,4) {$\vdots$};
		\node at (10,5) {$\circ$};

		\node at (11,1) {$\circ$};
		\node at (11,2) {$\circ$};
		\node at (11,3) {$\bullet$};
		\node at (11,4) {$\vdots$};
		\node at (11,5) {$\circ$};

	\end{tikzpicture}
		
		    \label{DomP3kplus1Km}
		\caption{Dominating Sets of $P_{3k+1} \times K_m$}
		\end{center}
		
        \end{figure}
		
		Let ~$B_0 = X_1 \cup X_2 $, \quad $B_i= \left(  X_{3i} \cup X_{3i+1} \cup X_{3i+2} \right)$ for $1 \leq i \leq k-1$,\\ and 	$B_k = X_{3k} \cup X_{3k+1}$. So that, $V =\left( \bigcup\limits_{i=0}^{k} B_i \right)$

		From $B_0$ we choose vertices $(1,3),(2,3)$.
		From each $B_i, ~ 1 \leq i \leq k-1$ we choose two vertices $(3i+1,j), (3i+2,j)$ where $
		j = \begin{cases}
			1 & \text{if $i$ is odd} \\
			2 & \text{if $i$ is even}
		\end{cases}.
		$ 
		
		From $B_k$ we choose the vertex $(3k, 3)(3k+1, 3)$.

		Using these vertices we obtain the set, $$D_{3k+1} =\{(1,3), (2,3), (4,1), (5,1), (7,2), (8,2), \ldots, (3k,3)(3k+1,3)\}.$$
		Observe that $D_{3k+1}$ is a dominating set for $P_{3k+1} \x K_m$.\\
		Hence, $\gamma(P_{3k+1} \x K_m) \leq 2k+2 $.

		\item Case 2. Let $n=3k+2$.\\ In this case we partition V as follows.
		
		$V =\left( \bigcup\limits_{i=0}^{k} B_i \right)$, where $B_0 = X_1 \cup X_2 $, and 
		
		$B_i=X_{3i} \cup X_{3i+1} \cup X_{3i+2}$ for $1 \leq i \leq k$.\\
		From $B_0$ we choose vertices $(1,3),(2,3)$.
		From each $B_i, ~ 1 \leq i \leq k-1$ we choose two vertices $(3i+1,j), (3i+2,j)$ where $
		j = \begin{cases}
			1 & \text{if $i$ is odd} \\
			2 & \text{if $i$ is even}
		\end{cases}.
		$\quad ~~\\ 
		From $B_k$ we choose vertices $(3k+1,3),(3k+2,3)$.\\~~

    \begin{figure}[H]

    \begin{center}
		\noindent 
%
%
%
	 
		\begin{tikzpicture}[scale=1]

			\node at (0,1) {$\tiny R_{1}$};
			\node at (0,2) {$R_{2}$};
			\node at (0,3) {$R_{3}$};
			\node at (0,4) {$\vdots$};
			\node at (0,5) {$R_{m}$};

			\node at (1,0) {$X_{1}$};
			\node at (2,0) {$X_{2}$};
			
			\draw (2.5,0.5) -- (2.5,5.5);
			
			\node at (3,0) {$X_{3}$};
			\node at (4,0) {$X_{4}$};
			\node at (5,0) {$X_{5}$};
			
			\draw (5.5,0.5) -- (5.5,5.5);

			\node at (6,0) {$X_{6}$};
			\node at (7,0) {$X_{7}$};
			\node at (8,0) {$X_{8}$};

			
			\draw (8.5,0.5) -- (8.5,5.5);

			\node at (10,0) {$X_{_{3k}}$};
			
			
			\node at (11,0) {$X_{_{3k+1}}$};
			
			%
			\node at (12,0) {$X_{_{3k+2}}$};

			\node at (1,1) {$\circ$};
			\node at (1,2) {$\circ$};
			\node at (1,3) {$\bullet$};
			\node at (1,4) {$\vdots$};
			\node at (1,5) {$\circ$};
			
			\node at (2,1) {$\circ$};
			\node at (2,2) {$\circ$};
			\node at (2,3) {$\bullet$};
			\node at (2,4) {$\vdots$};
			\node at (2,5) {$\circ$};
			
			\node at (3,1) {$\circ$};
			\node at (3,2) {$\circ$};
			\node at (3,3) {$\circ$};
			\node at (3,4) {$\vdots$};
			\node at (3,5) {$\circ$};

			\node at (4,1) {$\bullet$};
			\node at (4,2) {$\circ$};
			\node at (4,3) {$\circ$};
			\node at (4,4) {$\vdots$};
			\node at (4,5) {$\circ$};
			
			\node at (5,1) {$\bullet$};
			\node at (5,2) {$\circ$};
			\node at (5,3) {$\circ$};
			\node at (5,4) {$\vdots$};
			\node at (5,5) {$\circ$};

			\node at (6,1) {$\circ$};
			\node at (6,2) {$\circ$};
			\node at (6,3) {$\circ$};
			\node at (6,4) {$\vdots$};
			\node at (6,5) {$\circ$};

			\node at (7,1) {$\circ$};
			\node at (7,2) {$\bullet$};
			\node at (7,3) {$\circ$};
			\node at (7,4) {$\vdots$};
			\node at (7,5) {$\circ$};
			
			\node at (8,1) {$\circ$};
			\node at (8,2) {$\bullet$};
			\node at (8,3) {$\circ$};
			\node at (8,4) {$\vdots$};
			\node at (8,5) {$\circ$};
			
			
			
			\node at (9,1) {$\ldots$};
			\node at (9,5) {$\ldots$};
			\node at (9,4) {$\ddots$};
			
			\draw (9.5,0.5) -- (9.5,5.5);

			\node at (10,1) {$\circ$};
			\node at (10,2) {$\circ$};
			\node at (10,3) {$\circ$};
			\node at (10,4) {$\vdots$};
			\node at (10,5) {$\circ$};

			\node at (11,1) {$\circ$};
			\node at (11,2) {$\circ$};
			\node at (11,3) {$\bullet$};
			\node at (11,4) {$\vdots$};
			\node at (11,5) {$\circ$};
			
			\node at (12,1) {$\circ$};
			\node at (12,2) {$\circ$};
			\node at (12,3) {$\bullet$};
			\node at (12,4) {$\vdots$};
			\node at (12,5) {$\circ$};

		\end{tikzpicture}		
		
	     \caption{Dominating Sets of $P_{3k+2} \times K_m$}
        \end{center}
    \end{figure}
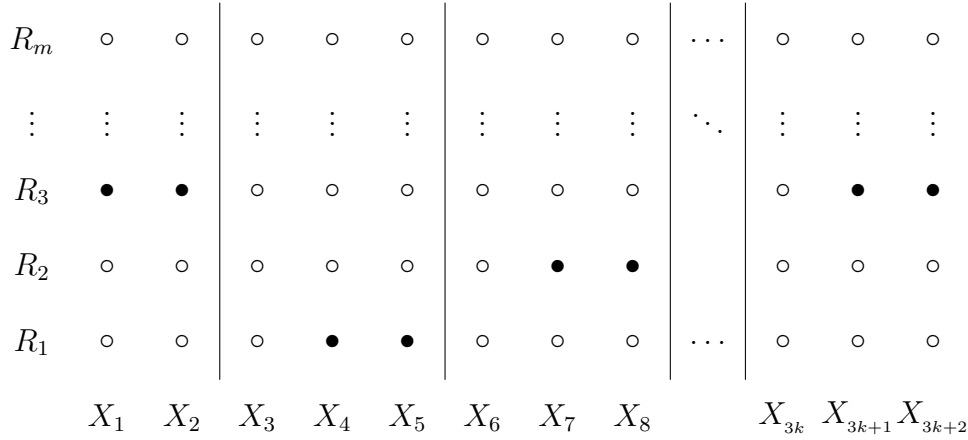	
		
		Using these vertices we obtain the set $$D_{3k+2} =\left\{(1,3), (2,3), (4,1), (5,1), (7,2), (8,2), \ldots, (3k+1, 3), (3k+2,3)\right\}$$
		Observe that $D_{3k+2}$ is a dominating set for $P_{3k+2} \x K_m$.\newline
		Hence, $\gamma(P_{3k+2} \x K_m) \leq 2k+2$.

		\item Case 3. Let $n=3k$.\\
		Similarly we obtain the set $D_{3k}$ as follows.
		
		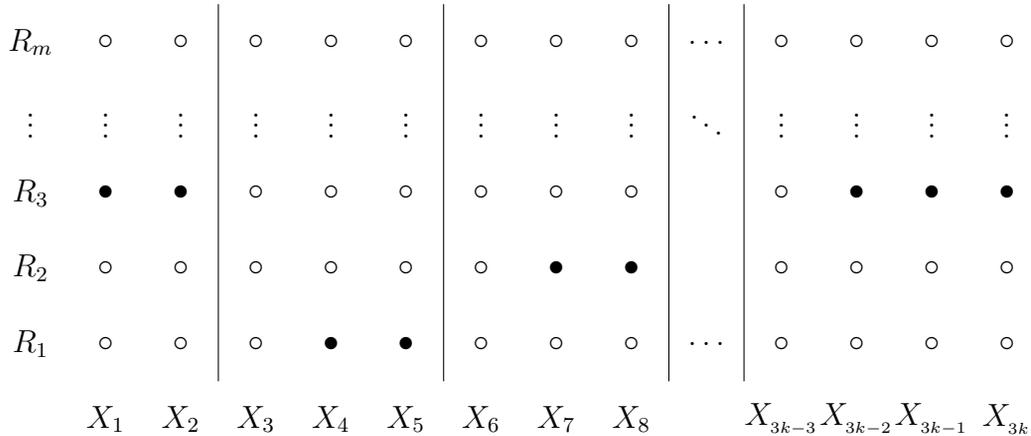
\begin{figure}[H]
		
		\begin{center}
			\noindent 
			
		\begin{tikzpicture}[scale=1]

			\node at (0,1) {$\tiny R_{1}$};
			\node at (0,2) {$R_{2}$};
			\node at (0,3) {$R_{3}$};
			\node at (0,4) {$\vdots$};
			\node at (0,5) {$R_{m}$};

			\node at (1,0) {$X_{1}$};
			\node at (2,0) {$X_{2}$};
			
			\draw (2.5,0.5) -- (2.5,5.5);
			
			\node at (3,0) {$X_{3}$};
			\node at (4,0) {$X_{4}$};
			\node at (5,0) {$X_{5}$};
			
			\draw (5.5,0.5) -- (5.5,5.5);

			\node at (6,0) {$X_{6}$};
			\node at (7,0) {$X_{7}$};
			\node at (8,0) {$X_{8}$};

			
			\draw (8.5,0.5) -- (8.5,5.5);

			\node at (10,0) {$X_{_{3k-3}}$};
			
			
			\node at (11,0) {$X_{_{3k-2}}$};
			
			%
			\node at (12,0) {$X_{_{3k-1}}$};
			
			\node at (13,0) {$X_{_{3k}}$};

			\node at (1,1) {$\circ$};
			\node at (1,2) {$\circ$};
			\node at (1,3) {$\bullet$};
			\node at (1,4) {$\vdots$};
			\node at (1,5) {$\circ$};
			
			\node at (2,1) {$\circ$};
			\node at (2,2) {$\circ$};
			\node at (2,3) {$\bullet$};
			\node at (2,4) {$\vdots$};
			\node at (2,5) {$\circ$};
			
			\node at (3,1) {$\circ$};
			\node at (3,2) {$\circ$};
			\node at (3,3) {$\circ$};
			\node at (3,4) {$\vdots$};
			\node at (3,5) {$\circ$};

			\node at (4,1) {$\bullet$};
			\node at (4,2) {$\circ$};
			\node at (4,3) {$\circ$};
			\node at (4,4) {$\vdots$};
			\node at (4,5) {$\circ$};
			
			\node at (5,1) {$\bullet$};
			\node at (5,2) {$\circ$};
			\node at (5,3) {$\circ$};
			\node at (5,4) {$\vdots$};
			\node at (5,5) {$\circ$};

			\node at (6,1) {$\circ$};
			\node at (6,2) {$\circ$};
			\node at (6,3) {$\circ$};
			\node at (6,4) {$\vdots$};
			\node at (6,5) {$\circ$};

			\node at (7,1) {$\circ$};
			\node at (7,2) {$\bullet$};
			\node at (7,3) {$\circ$};
			\node at (7,4) {$\vdots$};
			\node at (7,5) {$\circ$};
			
			\node at (8,1) {$\circ$};
			\node at (8,2) {$\bullet$};
			\node at (8,3) {$\circ$};
			\node at (8,4) {$\vdots$};
			\node at (8,5) {$\circ$};
			
			
			
			\node at (9,1) {$\ldots$};
			\node at (9,5) {$\ldots$};
			\node at (9,4) {$\ddots$};
			
			\draw (9.5,0.5) -- (9.5,5.5);

			\node at (10,1) {$\circ$};
			\node at (10,2) {$\circ$};
			\node at (10,3) {$\circ$};
			\node at (10,4) {$\vdots$};
			\node at (10,5) {$\circ$};

			\node at (11,1) {$\circ$};
			\node at (11,2) {$\circ$};
			\node at (11,3) {$\bullet$};
			\node at (11,4) {$\vdots$};
			\node at (11,5) {$\circ$};

		\node at (12,1) {$\circ$};
		\node at (12,2) {$\circ$};
		\node at (12,3) {$\bullet$};
		\node at (12,4) {$\vdots$};
		\node at (12,5) {$\circ$};

			\node at (13,1) {$\circ$};
			\node at (13,2) {$\circ$};
			\node at (13,3) {$\bullet$};
			\node at (13,4) {$\vdots$};
			\node at (13,5) {$\circ$};

		\end{tikzpicture}		
		
			\caption{Dominating Sets of $P_{3k} \times K_m$}
		    \label{DomP3kKm}
		\end{center}

        \end{figure}
    
		$$D_{3k} =\left\{(1,3), (2,3), (4,1), (5,1), (7,2), (8,2), \ldots, (3k-2,3),(3k-1,3),(3k,3)\right\}.$$
		Observe that $D_{3k}$ is a dominating set for $P_{3k} \x K_m$.\\
		Hence, $\gamma(P_{3k} \x K_m) \leq 2k+1$.	
	\end{itemize}	
	Therefore, the proposition is proved.
\end{proof}

%

The next theorem proves that the dominating sets obtained above are in fact the minimum dominating sets of $P_n \x K_m$.

\begin{theorem}\label{PLemma1} 
	For $m \geq 3,~n\geq 3$, 
	$$\gamma(P_n\times K_m) = \begin{cases} 
		2k+1 &\mbox{\textnormal{if}}~ n=3k,\\
		2k+2 &\mbox{\textnormal{if}}~ n=3k+1,~ 3k+2.
	\end{cases}$$ 
\end{theorem}	

\begin{proof}
	Let $D$ be a dominating set of $P_n \x K_m$.
	Using the fact, $d_0=d_{n+1}=0$ with Lemma \ref{lemma2}, we get the following results related to the corner pair of columns.
	
	\begin{enumerate}
		\item $d_1+d_2 \geq 2$.
		\item $d_{n-1}+d_{n} \geq 2$.
	\end{enumerate}

	Consider the following cases.
	\begin{itemize}

		\item Case 1. Let $n=3k+1$.\\
		Consider $$|D| = \sum_{i=1}^{3k+1} d_i = \sum_{i=0}^{3k+2} d_i = \sum_{i=0}^{k}\left(d_{3i}+d_{3i+1}+d_{3i+2}\right).$$
		Now by Lemma \ref{lemma2}, $\left(d_{3i}+d_{3i+1}+d_{3i+2}\right) \geq 2$, for $0 \leq i \leq k$. \\ Hence,  $|D| \geq 2(k+1).$

		\item Case 2. Let $n=3k+2$.\\
		Consider $$|D| = \sum_{i=1}^{3k+2} d_i = \sum_{i=0}^{3k+2} d_i = \sum_{i=0}^{k}\left(d_{3i}+d_{3i+1}+d_{3i+2}\right).$$
		Again by Lemma \ref{lemma2}, $\left(d_{3i}+d_{3i+1}+d_{3i+2}\right) \geq 2$, for $0 \leq i \leq k$.\\ Hence,  $|D| \geq 2(k+1).$
		
		\item Case 3. Let $n=3k$.\\
		We claim that $d_{3k}+d_{3k-1}+d_{3k-2}+d_{3k-3} \geq 3$.\\
		On the contrary, suppose $d_{3k-3}+d_{3k-2}+d_{3k-1}+d_{3k} < 3$.\\
		\emph{i.e.}  $d_{3k-3}+d_{3k-2}+d_{3k-1}+d_{3k} \leq 2$.\\
		But by Lemma \ref{lemma2}, $d_{3k-3}+d_{3k-2}+d_{3k-1} \geq 2$, which implies $d_{3k} = 0$.\\
		Similarly $d_{3k}+d_{3k-1} \geq 2$, implies $d_{3k-1}=d_{3k-2}=0$. \\
		Now column $X_{3k-1}$ is not dominated, which is a contradiction.\\
		Hence, $d_{3k}+d_{3k-1}+d_{3k-2}+d_{3k-3} \geq 3$.\\
		Again, $d_1+d_2 \geq 2$. \\
		Also, $d_{3i}+d_{3i+1}+d_{3i+2} \geq 2$, for $1 \leq i \leq k-2$.\\
		Adding all the inequalities above; we get \\
		$|D| \geq 2+2(k-2)+3$ \textit{i.e.}
		$|D| \geq 2k+1$.
	\end{itemize}
	Hence the theorem is proved using Proposition \ref{PLemma1}. \end{proof}
	
\subsection{Independent Domination Number and [1,2]-domination number of $C_n \x K_m$ and $P_n \x K_m$}	
	
	A subset of vertices is \textit{independent}  if no two vertices in it are adjacent. An \textit{independent dominating set} of $G$ is a set that is
	both dominating and independent in $G$. The \textit{independent domination number} $i(G)$ of $G$ is the minimum size of
	an independent dominating set. It follows immediately that $\gamma(G) \leq i(G)$. It is an open problem to characterize the graphs G such that $i(G) = \gamma(G)$.  For more details on such graphs see \cite{gupta2020graphs} and references therein. A subset of vertices is a \textit{$[1,2]$-set} if every vertex $v \in V$ is either in it or  adjacent to at least one but no more than two vertices in it. The \textit{$[1,2]$-domination number} of a graph $G$ $\gamma_{_{[1, 2]}} (G)$ is the minimum cardinality of a $[1, 2]$-set of $G$. Obviously, any [1,2]-set is also a dominating set, so $\gamma(G) \leq \gamma_{_{[1, 2]}} (G)$. Chellali \textit{et al.} \cite{chellali20131} proposed the following problem: for which graphs is $\gamma(G) = \gamma_{_{[1, 2]}} (G)$? Allan and Laskar \cite{allan1978domination} have shown that claw-free graphs are graphs with equal domination and independent domination numbers.
	Chellali \textit{et al.}. \cite{chellali20131} expanded on this by proving $\gamma(G) = i(G) = \gamma_{_{[1, 2]}} (G)$ for claw-free graph $G$.
	\par~\\ \noindent It is straightforward to verify that the dominating sets constructed by us are independent and [1,2]-dominating. Hence the following proposition is an immediate consequence.
	
	\begin{prop}\label{prop46}
		For $m \geq 2$, $n \geq 6$
		\begin{enumerate}
				\item $\gamma(C_n \x K_m) = i(C_n \x K_m) = \gamma_{_{[1, 2]}}(C_n \x K_m)$.
				\item $\gamma(P_n \x K_m) = i(P_n \x K_m) = \gamma_{_{[1, 2]}}(P_n \x K_m)$.
		\end{enumerate}	
	\end{prop} Hence, we get a large family of graphs with property domination number, independent domination number, and [1,2]-domination number being equal.

\subsection{Erratum}

In this section we address two mistakes that have appeared in the literature. T. Sitthiwirattham studied the domination number of direct product of path and cycle in \cite{Pnsitthiwirattham2012domination} and \cite{Cnsitthiwirattham2012domination} respectively. According to his claims, Theorem 2.3 in \cite{Pnsitthiwirattham2012domination} and Theorem 2.3 in  \cite{Cnsitthiwirattham2012domination}, for a graph $G$ of order $m$, 
$$\gamma(P_n \x G) = \gamma(C_n \x G) = \min ~ \left\{ n\gamma(G), m\left\lceil \frac{n}{3} \right\rceil    \right\}.$$
His result gives $\gamma(P_6 \x K_8) = 6 = \gamma(C_6 \x K_8)$. 
But $\gamma(P_6 \x K_8) = 5$ and $\gamma(C_6 \x K_8) = 4$. Also notice that, $\gamma(P_n \x K_m) $ need not be equal to $\gamma(C_n \x K_m)$. 
Hence Sitthiwirattham's claims are found to be invalid.\\

\par \noindent We also briefly mention an error in a paper  by Gravier and Khelladi \cite{gravier1995domination}. 
\begin{prop}[\cite{gravier1995domination}, Proposition 2.3]
	For $n > 1$ and every graph G we have $$\gamma(P_n \x G)  \leq 2 \gamma(G) \left( \left\lfloor \frac{n}{4} \right\rfloor +1\right)$$ 
\end{prop} \noindent For $n=6$, this gives us $\gamma(P_6 \x K_m) \leq 4$. But, we have proved that $\gamma(P_6 \x K_m) = 5$.\\  

\section{Secure Domination Number of \boldmath $C_n \x K_m$ and $P_n \x K_m$ }
\subsection{Secure Domination Number of \boldmath $C_n \x K_m$ }

For a secure dominating set $S$ of either $P_n \x K_m$ or  $C_n \x K_m$, and for $i \in [n]$, let $s_i$ denote the cardinality of $S \cap X_i$. As mentioned before, in $P_n \x K_m$, it is convenient to assume that  $s_0=s_{n+1}=0$, whereas in case of $C_n \x K_m$, we suppose $s_{0}=s_{n}$ and $s_{n+1}=s_{1}$.\\

\noindent The next lemma is crucial in the proof of subsequent theorems.
\begin{lemma} \label{Slemma2}
	Let $m \geq 3, n \geq 2$. If $S$ be a secure dominating set of either $C_n \x K_m$ or $P_n \x K_m$, then  $s_{i-1}+s_{i}+s_{i+1} \geq 3$, for $1 \leq i \leq n$. 
\end{lemma}
\begin{proof}
	By Lemma \ref{lemma2},  $s_{i-1}+s_{i}+s_{i+1} \geq 2$. On the contrary, if possible, let $s_{i-1}+s_{i}+s_{i+1} =2$. This equation has six possible solutions. Let us analyze them one by one.
	\begin{enumerate}[(i)]
		
		\item Let $s_{i-1}=0$, $s_i=2$, $s_{i+1}=0$.\\
		This is not a dominating set since it 
		violets condition (a) of Lemma \ref{SLemma1}.

		\item Let $s_{i-1}=0$, $s_i=0$, $s_{i+1}=2.$\\
		In this case,  swap any dominating vertex in $X_{i+1}$ with a vertex in $X_i$, which is not in the same row as that of vertices in $X_{i+1}$. The new set after replacement is not a dominating set by Lemma \ref{SLemma1}.

\begin{figure}[H]

			\begin{center}
			\noindent 
			\begin{tikzpicture}
		
		\node at (0,3) {$\circ$};
		\node at (1,3) {$\bullet$};
		\node at (2,3) {$\circ$};

		\node at (-2,3) {$\circ$};
		\node at (-3,3) {$\circ$};
		\node at (-4,3) {$\circ$};

		\node at (-2,2) {$\vdots$};
		\node at (-3,2) {$\vdots$};
		\node at (-4,2) {$\vdots$};
		
		\node at (-2,1) {$\bullet$};
		\node at (-3,1) {$\circ$};
		\node at (-4,1) {$\circ$};
		
		\node at (-2,0) {$\bullet$};
		\node at (-3,0) {$\circ$};
		\node at (-4,0) {$\circ$};
		
		\node at (-4,-1)  {Initial dominating set of $X_i$};
		
		\node at (-2,-0.5)  {$X_{i+1}$};
		\node at (-3,-0.5)  {$X_{i}$};
		\node at (-4,-0.5)  {$X_{i-1}$};
		
		\draw (-1,-0.5) -- (-1,3.5);
		
		\node at (0,3) {$\circ$};
		\node at (1,3) {$\circ$};
		\node at (2,3) {$\circ$};

		\node at (0,2) {$\vdots$};
		\node at (1,2) {$\vdots$};
		\node at (2,2) {$\vdots$};

		\node at (0,1) {$\circ$};
		\node at (1,1) {$\circ$};
		\node at (2,1) {$\bullet$};
		
		\node at (0,0) {$\circ$};
		\node at (1,0) {$\circ$};
		\node at (2,0) {$\circ$};
		
		\node at (0,-0.5) {$X_{i-1}$};
		\node at (1,-0.5) {$X_{i}$};
		\node at (2,-0.5) {$X_{i+1}$};
		
		\node at (3,-1)  {New sets after swapping};

		\node at (4,3) {$\circ$};
		\node at (5,3) {$\bullet$};
		\node at (6,3) {$\circ$};

		\node at (4,2) {$\vdots$};
		\node at (5,2) {$\vdots$};
		\node at (6,2) {$\vdots$};

		\node at (4,1) {$\circ$};
		\node at (5,1) {$\circ$};
		\node at (6,1) {$\circ$};
		
		\node at (4,0) {$\circ$};
		\node at (5,0) {$\circ$};
		\node at (6,0) {$\bullet$};
		
		\node at (4,-0.5) {$X_{i-1}$};
		\node at (5,-0.5) {$X_{i}$};
		\node at (6,-0.5) {$X_{i+1}$};
		
	\end{tikzpicture}		
	
             \caption{Figure Corresponding to Case (ii)}
    \label{s1s2s3case1}

			\end{center}
		
\end{figure}
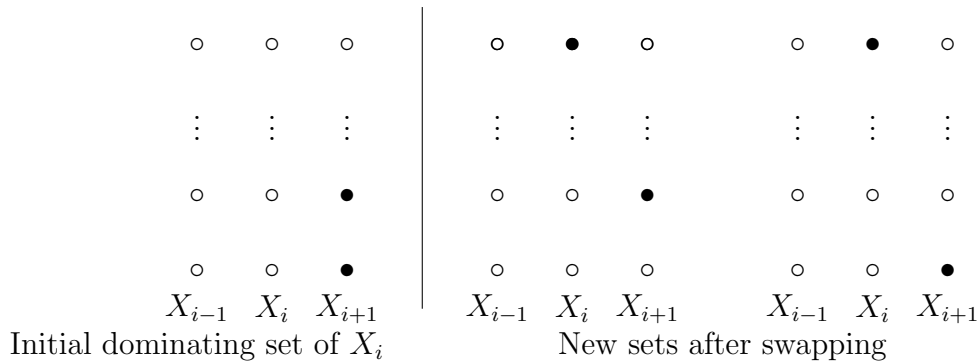
		\item $s_{i-1}=2$, $s_i=0$, $s_{i+1}=0$.\\
		This case is similar case (ii) above, obtained replacing $(i+1)$ by $(i-1)$.
		
		\item $s_{i-1}=0$, $s_i=1$, $s_{i+1}=1$.\\
		Swapping the dominating vertex in $X_{i+1}$ with any vertex in $X_i$, produces a new set, which is a subset of $X_i$. The new set after replacement is not a dominating since it violets condition (a) of Lemma \ref{SLemma1}.

  \item $s_{i-1}=1$, $s_i=1$, $s_{i+1}=0$.\\
 This case is similar case (iv) above, obtained replacing $(i+1)$ by $(i-1)$.

        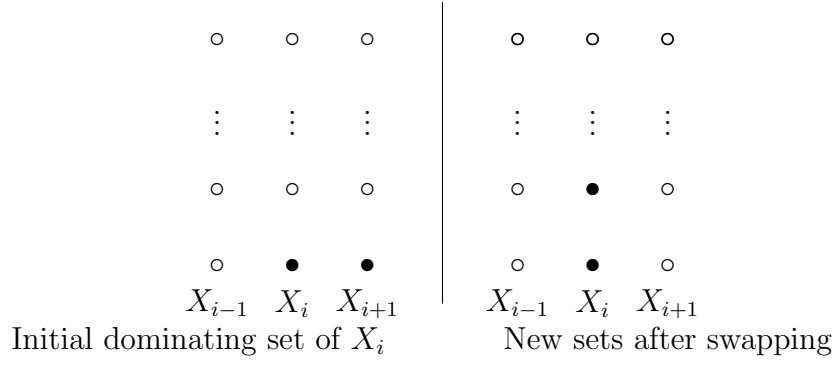
\begin{figure}[H]
        \begin{center}
			\noindent 
				\begin{tikzpicture}
			
			\node at (0,3) {$\circ$};
			\node at (1,3) {$\circ$};
			\node at (2,3) {$\circ$};

			\node at (-2,3) {$\circ$};
			\node at (-3,3) {$\circ$};
			\node at (-4,3) {$\circ$};

			\node at (-2,2) {$\vdots$};
			\node at (-3,2) {$\vdots$};
			\node at (-4,2) {$\vdots$};
			
			\node at (-2,1) {$\circ$};
			\node at (-3,1) {$\circ$};
			\node at (-4,1) {$\circ$};
			
			\node at (-2,0) {$\bullet$};
			\node at (-3,0) {$\bullet$};
			\node at (-4,0) {$\circ$};
			
			\node at (-4.25,-1)  {Initial dominating set of $X_i$};
			
			\node at (-2,-0.5)  {$X_{i+1}$};
			\node at (-3,-0.5)  {$X_{i}$};
			\node at (-4,-0.5)  {$X_{i-1}$};
			
			\draw (-1,-0.5) -- (-1,3.5);
			
			\node at (0,3) {$\circ$};
			\node at (1,3) {$\circ$};
			\node at (2,3) {$\circ$};

			\node at (0,2) {$\vdots$};
			\node at (1,2) {$\vdots$};
			\node at (2,2) {$\vdots$};

			\node at (0,1) {$\circ$};
			\node at (1,1) {$\bullet$};
			\node at (2,1) {$\circ$};
			
			\node at (0,0) {$\circ$};
			\node at (1,0) {$\bullet$};
			\node at (2,0) {$\circ$};
			
			\node at (0,-0.5) {$X_{i-1}$};
			\node at (1,-0.5) {$X_{i}$};
			\node at (2,-0.5) {$X_{i+1}$};
			
			\node at (2,-1)  {New sets after swapping};

%
%
%
%
%
%
			
		\end{tikzpicture}		
	
            \caption{Figure Corresponding to Case (iv)}
    \label{s1s2s3case4}
		\end{center}
		\end{figure}
	  
		\item $s_{i-1}=1$, $s_i=0$, $s_{i+1}=1$.\\
		Swap any dominating vertex with a vertex $X_i$, which does not lie in same row as that of dominating vertices. The resulting set is not dominating by 
		Lemma \ref{SLemma1}.

		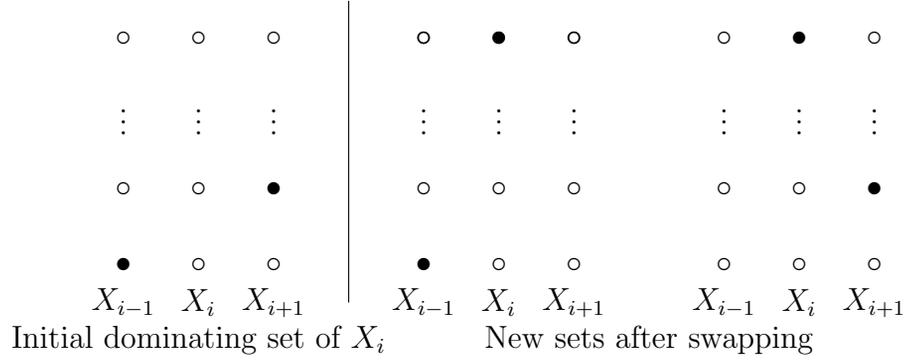
\begin{figure}[H]
         
			\begin{center}
			\noindent 
			
		\begin{tikzpicture}
			
			\node at (0,3) {$\circ$};
			\node at (1,3) {$\circ$};
			\node at (2,3) {$\circ$};

			\node at (-2,3) {$\circ$};
			\node at (-3,3) {$\circ$};
			\node at (-4,3) {$\circ$};

			\node at (-2,2) {$\vdots$};
			\node at (-3,2) {$\vdots$};
			\node at (-4,2) {$\vdots$};
			
			\node at (-2,1) {$\bullet$};
			\node at (-3,1) {$\circ$};
			\node at (-4,1) {$\circ$};
			
			\node at (-2,0) {$\circ$};
			\node at (-3,0) {$\circ$};
			\node at (-4,0) {$\bullet$};
			
			\node at (-3,-1)  {Initial dominating set of $X_i$};
			
			\node at (-2,-0.5)  {$X_{i+1}$};
			\node at (-3,-0.5)  {$X_{i}$};
			\node at (-4,-0.5)  {$X_{i-1}$};
			
			\draw (-1,-0.5) -- (-1,3.5);
			
			\node at (0,3) {$\circ$};
			\node at (1,3) {$\bullet$};
			\node at (2,3) {$\circ$};

			\node at (0,2) {$\vdots$};
			\node at (1,2) {$\vdots$};
			\node at (2,2) {$\vdots$};

			\node at (0,1) {$\circ$};
			\node at (1,1) {$\circ$};
			\node at (2,1) {$\circ$};
			
			\node at (0,0) {$\bullet$};
			\node at (1,0) {$\circ$};
			\node at (2,0) {$\circ$};
			
			\node at (0,-0.5) {$X_{i-1}$};
			\node at (1,-0.5) {$X_{i}$};
			\node at (2,-0.5) {$X_{i+1}$};
			
			\node at (3,-1)  {New sets after swapping};

			\node at (4,3) {$\circ$};
			\node at (5,3) {$\bullet$};
			\node at (6,3) {$\circ$};

			\node at (4,2) {$\vdots$};
			\node at (5,2) {$\vdots$};
			\node at (6,2) {$\vdots$};

			\node at (4,1) {$\circ$};
			\node at (5,1) {$\circ$};
			\node at (6,1) {$\bullet$};
			
			\node at (4,0) {$\circ$};
			\node at (5,0) {$\circ$};
			\node at (6,0) {$\circ$};
			
			\node at (4,-0.5) {$X_{i-1}$};
			\node at (5,-0.5) {$X_{i}$};
			\node at (6,-0.5) {$X_{i+1}$};
			
		\end{tikzpicture}		
            \caption{Figure Corresponding to Case (vi)}
    \label{s1s2s3case5}
		\end{center}
		\end{figure}

	\end{enumerate}
	Thus, if $S$ is a secure dominating subset of $V$ then $s_{i-1}+s_{i}+s_{i+1} =2$ is not possible.  Hence, $s_{i-1}+s_{i}+s_{i+1} \geq 3$, for $1 \leq i \leq n$.
\end{proof}

We now state the main theorem of this section that determines the exact value of secure domination number of 	$C_n\times K_m$.

\begin{theorem} 
	For $m \geq 3$, 
	(i) $\gamma_s(C_n\times K_m) =  n$, \textnormal{if} $n \geq 3$.\\
	(ii) $\gamma_s(C_2\times K_m) =  \begin{cases} 
		3 &\mbox{\textnormal{if}}~ m=3,\\
		4 &\mbox{\textnormal{if}}~ m \geq 4.
	\end{cases}$
\end{theorem}

\begin{proof} (i) Let $m \geq 3, n\geq 3$. Consider the following cases.
	\begin{itemize}
		\item Case 1: Let $n=3k$. By lemma \ref{Slemma2}, $s_{3i-2}+s_{3i-1}+s_{3i} \ge 3$ for $i=1,2, \cdots k$. So $|S| \geq 3k$
		\item Case 2: Let $n=3k+1$.\\ Let, if possible, assume that $|S| \leq 3k$. \\
		 Now, $$|S| = \sum_{i=1}^{3k+1} \left(s_{3i-2}+s_{3i-1}+s_{3i}\right) \leq 3k.$$  Choosing the three consecutive columns partition of $C_{3k+1} \times K_m$ that excludes $X_j$,  by lemma \ref{Slemma2} we get, $$\left(\sum_{i=1}^{3k+1} s_i\right) - s_j \geq 3k.$$ 
		 The two inequalities together imply $s_j=0$. Now, since $j$ is arbitrary, we arrive at a contradiction. So, $|S| \geq 3k+1$
		\item Case 3: Let $n=3k+2$.\\Let, if possible, assume that, $S$ is a secure dominating set such that $|S| \leq 3k+1$.
		Since $n=3k+2$, excluding   columns  $X_{t+1}$ and $X_{t+2}$, group remaining $3k$ columns in triples of three consecutive columns. 
		Now, $\textnormal{for each}~~ t \in [3k+1]$, we have
		$$|S| = \sum_{i=1}^{3k+2} s_i = \left(\sum_{\substack{i=1 \\i \neq t+1, ~t+2}}^{3k+2} s_i \right) + \left(s_{t+1} + s_{t+2}\right)= \left(\sum_{\substack{i=1 \\i \neq t+1, ~t+2}}^{3k+2} s_i\right) + \left(s_{t+1} + s_{t+2} +s_t\right) - s_t,$$
		\noindent By Lemma \ref{Slemma2}, $ |S| \geq 3k + 3 -s_t$ ~~ $\textnormal{for each~} t \in [3k+1].$ \\In particular, since we assumed that	$|S| \leq 3k+1$, there exists a $j \in [3k+1]$ such that $s_j = 0$.This implies $ |S| \geq 3k+3,$ which is a contradiction to the initial assumption $|S| \leq 3k+1$.  So, $|S| \geq 3k+2$.
	\end{itemize}
	
	 \noindent Therefore, in any case,  $\gamma_s(C_n\times K_m)\geq  n.$ \\~\par Now, It is easy to see that the first row viz. $R_1$ is a secure dominating set.\\ So, $\gamma_s(C_n\times K_m)\leq  n.$ \\~~
	
	\par \noindent Hence, for $m \geq 3, n \geq 3$ the secure domination  number, $\gamma_s(C_n\times K_m) =  n.$ \\
	\par
	(ii) Let $m \geq 3, n=2$. Now by Lemma \ref{Slemma2}, $|S| =s_1 + s_2 \geq 3$.
	If $m = 3$ then $X_1$ is a dominating set. So, $\gamma_s(C_2 \x K_3) = 3$.
	Now suppose,  $m \geq 4$.
	If $s_1=0$ and $s_2=3$, then the fourth vertex in $X_1$ is not dominated.
	If $s_1=1$ and $s_2=2$, then after swapping dominating vertex in $X_1$, the new set lies in $X_2$ which is not dominating as before.
	The other two cases are similar. Hence there does not exists a secure dominating set of cardinality 3.
	Now it is straighforward to verify that $R_1 \cup R_2$ is a secure dominating set of $C_2 \x K_m$. So $\gamma_s(C_2 \x K_m) = 4$, if $m \geq 4$.
	\par \noindent Hence, $\gamma_s(C_2\times K_m) =  \begin{cases} 
		3 &\mbox{\textnormal{if}}~ m=3,\\
		4 &\mbox{\textnormal{if}}~ m \geq 4.
	\end{cases}$
\end{proof}

\subsection{Secure Domination Number of \boldmath $P_n \x K_m$}

In this section we study the secure dominating sets of $P_n \x K_m$.

\begin{lemma}\label{Corner31}  For $m \geq 4,~n \geq 3$, if $S$ is a secure dominating set of $P_n \x K_m$, then \\(i) $s_1+ s_2 + s_3 \geq 4$.\\
(ii)  $s_n+ s_{n-1} + s_{n-2} \geq 4$.
\end{lemma}
\begin{proof}
	(i)Since in $P_n \x K_m$, $s_0 = 0$, $s_1+ s_2 \geq 3$, by Lemma \ref{Slemma2}. Also, $s_1+ s_2 + s_3 \geq 3$, by Lemma \ref{Slemma2}. Let, if possible, $s_1+ s_2 + s_3 = 3$. So $s_1+ s_2 = 3$ and $s_3 =0$. The equation $s_1+ s_2 = 3$ has four possible solutions. 
	\begin{itemize}
		\item Case 1: $s_1=3, s_2=0$. In this case the fourth vertex in $X_1$, other than the three in $S$, is not dominated by $S$.
		\item Case 2: $s_1=0, s_2=3$. In this case the fourth vertex in $X_2$, other than the three in $S$, is not dominated by $S$.
		\item Case 3: $s_1=1, s_2=2$. In this case, swap any vertex in $X_2$, (other than two dominating vertices),  with dominating vertex in $X_1$. The swap set does not dominate the fourth vertex in $X_2$, (other than the three in swap set).
		\item Case 4: $s_1=2, s_2=1$. In this case, swap any vertex in $X_1$, (other than two dominating vertices),  with dominating vertex in $X_2$. The swap set does not dominate the fourth vertex in $X_1$, (other than the three in swap set).
	\end{itemize}
	So, $s_1+ s_2 + s_3 = 3$ is not possible. Hence, $s_1+ s_2 + s_3 \geq 4.$ \\
	\par (ii) The proof is same as above.	
\end{proof}

The lemmas \ref{Slemma2} and \ref{Corner31} enables us to establish the main result of this section which determines the exact value of secure domination number of 	$P_n\times K_m$.

\begin{theorem} Let $n \geq 3$. (i) $\gamma_s(P_n\times K_m) =  n+2$ \textnormal{if} $m \geq 4$. \\
	(ii) $\gamma_s(P_n\times K_3) = \begin{cases} 
		3k &\mbox{\textnormal{if}}~ n=3k,\\
		3k+3 &\mbox{\textnormal{if}}~ n=3k+1,\\
		3k+3 &\mbox{\textnormal{if}}~ n=3k+2.\\
	\end{cases}$
	
\end{theorem}
\begin{proof} (i) Let $m \geq 4$. Consider the following cases.
	\begin{itemize}
		\item Case 1:  Let $n=3k$.\\
		By Lemma \ref{Slemma2}, $$\sum_{i=2}^{k-1} (s_{3i-2}+s_{3i-1+}+s_{3i})  \geq 3(k-2).$$
		By Lemma \ref{Corner31}, $s_1+ s_2 + s_3 \geq 4$ and	$s_{3k}+ s_{3k-1} + s_{3k-2} \geq 4$.\\
		Adding three inequalities we get $|S| \geq 3k+2$.
		
		\item Case 2: Let $n=3k+1$.\\
        By Lemma \ref{Slemma2}, $$\sum_{i=1}^{k-1} (s_{3i}+s_{3i+1+}+s_{3i+2})  \geq 3(k-2).$$
		By Lemma \ref{Corner31}, $s_1+ s_2  \geq 3$ 
		and	$s_{3k+1}+ s_{3k} \geq 3$.\\
		Adding three inequalities we get $|S| \geq 3k+3$.
		
		\item Case 3: Let $n=3k+2$.\\
	    By Lemma \ref{Slemma2}, $$\sum_{i=1}^{k-1} (s_{3i}+s_{3i+1+}+s_{3i+2})  \geq 3(k-2).$$
	    By Lemma \ref{Corner31}, $s_1+ s_2  \geq 3$ 
	    and	$s_{3k+2}+ s_{3k+1} + s_{3k} \geq 4$.\\
		Adding three inequalities we get $|S| \geq 3k+4$.
	\end{itemize}
	\par \noindent Therefore, in any case, $\gamma_s(P_n\times K_m)\geq  n+2.$ \\
	\par Now it is easy to see that $R_1 \cup \{((2,2), (n-1,2))\}$ is a secure dominating set. \\
	So, $\gamma_s(P_n\times K_m)\leq  n+2.$~~ Hence, $\gamma_s(P_n\times K_m) =  n+2,$ whenever $m \geq 4$. \\~~\\~
	
	\par (ii) 
Consider the graph $P_3 \x K_3$. By Lemma \ref{Slemma2}, $\gamma_s(P_3 \x K_3) \geq 3$. Now the middle column $X_2$ is clearly a secure dominating set of $P_3 \x K_3$, gives us $\gamma_s(P_3 \x K_3) = 3$. We replicate this middle column in block of three columns to get a dominating set.  

	\begin{itemize}
	\item Case 1: Let $n=3k$.\\
	We divide the vertex set into blocks of three columns each we get, and choosing the middle column from each block, we get a secure dominating set, Thus $\gamma_s(P_{3k} \x K_3) \leq 3k$. But By Lemma \ref{Slemma2}, $$|S| = \sum_{i=1}^{k} (s_{3i-2}+s_{3i-1+}+s_{3i})  \geq 3k.$$ Thus $\gamma_s(P_{3k} \x K_3) = 3k$.
    
	\item Case 2: Let $n=3k+1$.\\
	Except for the last column, we divide the vertex set into blocks of three columns each we get, and choosing the middle column from each block along with the last column, we get a secure dominating set, Thus $\gamma_s(P_{3k+1} \x K_3) \leq 3k+1$.\\
	Now, as $s_0=s_{3k+2}=0$, by Lemma \ref{Slemma2} $$|S| =  \sum_{i=0}^{k} (s_{3i}+s_{3i+1+}+s_{3i+2}) \geq 3(k+1).$$
	Thus $\gamma_s(P_{3k+1} \x K_3) = 3k+3$.
	
	\item Case 3: Let $n=3k+2$.\\
	Except for the last two columns, we divide the vertex set into blocks of three columns each we get, and choosing the middle column from each block along with the last column, we get a secure dominating set, Thus $\gamma_s(P_{3k+2} \x K_3) \leq 3k+3$. Now as $s_{3k+3}=0$ by Lemma \ref{Slemma2}, $$|S| = \sum_{i=1}^{k+1} (s_{3i-2}+s_{3i-1+}+s_{3i})  \geq 3(k+1).$$ Thus $\gamma_s(P_{3k+2} \x K_3) = 3k+3$.
\end{itemize}

\noindent Hence the theorem is proved. \end{proof} 

\subsection{2-Domination Number of \boldmath $P_n \x K_m$ and $C_n \x K_m$   }
A vertex subset of a graph $G$ is said to 2-dominate the graph if each vertex of $G$ is either in the given subset or  has
at least two neighbors in it. The minimum cardinality of a
2-dominating set, denoted by $\gamma_2(G)$,  is called the 2-domination number of the graph $G$.
A 2-dominating set is clearly a secure dominating set and so, $\gamma_2(G) \geq   \gamma_s(G)$. It is a natural to ask: for which graphs is $\gamma_2(G) = \gamma_s(G)$?

	\par~\\ \noindent It is straightforward to verify that the secure dominating sets constructed by us are, in fact, 2-dominating. Hence the following proposition is an immediate consequence.

\begin{prop}
	For $m \geq 3$, $n \geq 2$
	\begin{enumerate}
		\item $\gamma_s(C_n \x K_m) = \gamma_2(C_n \x K_m)$.
		\item $\gamma_s(P_n \x K_m) = \gamma_2(P_n \x K_m)$.
	\end{enumerate}	
\end{prop} 
Hence, we get a family of graphs with equal  2-domination number and secure domination number.

\section{Conclusion}

In this paper, we obtained the exact value of the  domination number  of the direct product of a path or a cycle with a complete graph. The determination of the family of graphs $G$ for which $i(G) = \gamma(G)$ remains an open problem, as discussed in \cite{allan1978domination}.
Chellali et al. \cite{chellali20131} posed the question of identifying a family of graphs where $\gamma(G) = \gamma_{_{[1, 2]}} (G)$. Allan and Laskar \cite{allan1978domination} demonstrated that claw-free graphs exhibit equal domination and independent domination numbers. Chellali et al. \cite{chellali20131} extended this result by proving $\gamma(G) = i(G) = \gamma_{_{[1, 2]}} (G)$ for claw-free graphs $G$.  Chellali et al. \cite{chellali20131} asked to obtain family of graphs for which graphs is $\gamma(G) = \gamma_{_{[1, 2]}} (G)$. Allan and Laskar \cite{allan1978domination} have shown that claw-free graphs are graphs with equal domination and independent domination numbers.	Chellali et al. \cite{chellali20131} expanded on this by proving $\gamma(G) = i(G) = \gamma_{_{[1, 2]}} (G)$ for claw-free graph $G$.   For the family of graphs onsisting of  direct product of a path or a cycle with a complete graph, it is shown that the independent domination number and [1,2]-domination number coincide with the domination number. We thus obtained a partial solution to these problems, distinct from the previously explored case of claw-free graphs \cite{chellali20131}. As a consequence of our work, counterexamples are provided to disprove some erroneous results obtained by T. Sitthiwirattham  in \cite{Pnsitthiwirattham2012domination} and \cite{Cnsitthiwirattham2012domination}. Furthermore, we determined the exact values of the secure domination number and the 2-domination number for the aforementioned graph family, establishing their equality. Hence, we identified a family of graphs with equal 2-domination number and secure domination number.


\section*{Acknowledgement}
The first author gratefully acknowledges the Council of Scientific and Industrial Research, New Delhi, India for the award of the Senior Research Fellowship (09/137(0617)/2019-EMR-I). The second/corresponding author acknowledges the Department of Science and Technology, New Delhi, India for the award of Women Scientist Scheme (DST/WOS-A/PM-14/2021(G)) for reseach in Basic/Applied Sciences.

\bibliographystyle{plain}
\bibliography{references}	

\end{document}